\documentclass{book}
\usepackage[contrib,lang=british,counterwithoutchapter]{ems-book-no-chap} 

\usepackage[textsize = tiny]{todonotes}
\usepackage[capitalize]{cleveref}
\usepackage[all, cmtip]{xy}
\usepackage{makecell}
\usepackage{verbatim}

\newtheorem{theorem}{Theorem}

\theoremstyle{definition}

\newtheorem{conjecture}[theorem]{Conjecture}

	\newcommand{\ls}{\left\{}
	\newcommand{\rs}{\right\}}
	\newcommand{\sm}{\setminus}
	
	\newcommand{\mcA}{\ensuremath{\mathcal{A}}}

	\newcommand{\mcM}{\ensuremath{\mathcal{M}}}
	
	\newcommand{\mcO}{\ensuremath{\mathcal{O}}}

	\newcommand{\mbQ}{\ensuremath{\mathbb{Q}}}
	
	\newcommand{\mbZ}{\ensuremath{\mathbb{Z}}}
	\newcommand{\mbN}{\ensuremath{\mathbb{N}}}
	\newcommand{\mbK}{\ensuremath{\mathbb{K}}}

	\newcommand{\on}[1]{\operatorname{#1}}

	\newcommand{\tildee}[1]{\widetilde{#1}}
	
	\newcommand{\Aut}{\ensuremath{\operatorname{Aut}}}
	\newcommand{\Out}{\ensuremath{\operatorname{Out}}}
	
	\newcommand{\Sym}{\operatorname{Sym}}
	\newcommand{\St}{\operatorname{St}}
	
	\newcommand{\MCG}{\ensuremath{\operatorname{MCG}}}

	\newcommand{\rk}{\operatorname{rk}}

	\newcommand{\vcd}{\operatorname{vcd}}

	\newcommand{\GL}[2]{\ensuremath{\operatorname{GL}_{#1}(#2)}}
	\newcommand{\SL}[2]{\ensuremath{\operatorname{SL}_{#1}(#2)}}
	\newcommand{\Sp}[2]{\ensuremath{\operatorname{Sp}_{#1}(#2)}}

    \newcommand{\tA}{\ensuremath{\mathtt{A}}}
    \newcommand{\tB}{\ensuremath{\mathtt{B}}}
    \newcommand{\tC}{\ensuremath{\mathtt{C}}}
    \newcommand{\tD}{\ensuremath{\mathtt{D}}}
    \newcommand{\tE}{\ensuremath{\mathtt{E}}}
    \newcommand{\tF}{\ensuremath{\mathtt{F}_4}}
    \newcommand{\tG}{\ensuremath{\mathtt{G}_2}}

\NewDocumentCommand {\Irel} { O{n} O{m} }{\mathcal{I}^{#2}_{#1}}

\NewDocumentCommand {\Idelrel} { O{n} O{m} }{\mathcal{I}^{\delta,#2}_{#1}}

\NewDocumentCommand {\Isigdelrel} { O{n} O{m} }{\mathcal{I}^{\sigma,\delta,#2}_{#1}}

\NewDocumentCommand {\IArel} { O{n} O{m} }{\mathcal{IA}^{#2}_{#1}}

\NewDocumentCommand {\IAArel} { O{n} O{m} }{\mathcal{IAA}^{#2}_{#1}}

\NewDocumentCommand {\IAAstrel} { O{n} O{m} }{\mathcal{IAA}^{*,#2}_{#1}}

    \newcommand{\field}{\ensuremath{\mbK}}
	\newcommand{\ring}{\ensuremath{R}}
	\newcommand{\numberring}{\mcO}
	\newcommand{\chevalley}{\ensuremath{\mathcal{G}}} 
\newcommand{\class}{\operatorname{cl}}
\newcommand{\group}{\ensuremath{\Gamma}}

\newcommand{\building}{\Delta}

\newcommand{\CoxeterGroupTypeB}{\on{Sym}^\pm}

\begin{document}
\mainmatter

\title{(Non-)Vanishing of high-dimensional group cohomology}
\titlemark{High-dimensional group cohomology}

\emsauthor{1}{Benjamin Br{\"u}ck}{B.~Br{\"u}ck}


\emsaffil{1}{Institut f{\"u}r Mathematische Logik und Grundlagenforschung, Einsteinstr. 62
48149 M{\=u}nster, Germany \email{benjamin.brueck@uni-muenster.de}}

\classification[11F75, 20F28, 57K20]{57M07}

\keywords{group cohomology, duality groups, Steinberg module}


\begin{abstract}
Church--Farb--Putman formulated stability and vanishing conjectures for the high-dimensional cohomology of $\SL{n}{\mbZ}$, surface mapping class groups and automorphism groups of free groups. This is a survey on the current status of these conjectures and their generalisations.
\end{abstract}

\makecontribtitle


\section{Introduction}
This article is concerned with the rational cohomology of groups ``similar to $\SL{n}{\mbZ}$'' such as the surface mapping class group $\MCG(\Sigma_g)$ and the automorphism group of the free group $\Aut(F_n)$.
Each such group $\group$ has finite virtual cohomological dimension $\vcd(\group)\in \mbN$. This implies that its rational cohomology can only be non-trivial in finitely many degrees because it vanishes in degrees above $\vcd(\group)$,
\begin{equation*}
	H^k(\group;\mbQ) = 0 \text{ for } k > \vcd(\group).
\end{equation*}
Homological stability results give a comparably good understanding of the low-dimensional cohomology ($k\ll \vcd(\group)$) of these groups.
However, much less is known about their high-dimensional co\-ho\-mo\-logy ($k \approx \vcd(\group)$).
In \cite{CFP:stabilityconjectureunstable}, Church--Farb--Putman conjectured certain patterns in this unstable cohomology.

\paragraph{Duality}
The point of departure for these conjectures is that $\SL{n}{\mbZ}$, $\MCG(\Sigma_g)$ and $\Aut(F_n)$ are all virtual Bieri--Eckmann duality groups, which means that they satisfy an analogue of Poincar\'e duality. 
For $\SL{n}{\mbZ}$, this
follows from work by Borel--Serre \cite{BS:Cornersarithmeticgroups}: Let $\ring$ be a number ring, i.e.~the ring of integers in a number field $\field$, and $\chevalley$ a Chevalley--Demazure group scheme (for example $\chevalley = \on{SL}_n$, $\on{Sp}_{2n}$, $\on{PGL}_n$, $\on{Spin}_{2n+1}$ or $\on{SO}_{2n}$). Then $\group = \chevalley(\ring)$ is an arithmetic subgroup of $\chevalley(\field)$ and Borel--Serre's result implies that for all codimensions $0\leq i\leq \vcd(\group)<\infty$, there is an isomorphism
\begin{equation}
\label{eq:BE_duality}
H^{\vcd(\group)-i}(\group; \mbQ) \cong H_i(\group;\St(\Gamma) \otimes \mbQ).
\end{equation}
Here $\St(\group)$ is the \emph{Steinberg module}, i.e.~the top-degree
homology of the Tits building\footnote{This building is a highly symmetric simplicial $\Gamma$-complex. For $\Gamma = \SL{n}{\ring}$ and $\Sp{2n}{\ring}$, definitions are given in \cref{sec_resolutions_form} and \cref{sec_res_symplectic}.} $\building(\group)$ associated to $\chevalley(\field)$.
The duality given by \cref{eq:BE_duality} allows one to study high-dimensional cohomology with rational coefficients by investigating low-dimensional homology with coefficients in the dualising module 
$\St(\group)$. 
The trade-off is that the latter requires a good understanding of $\St(\group)$, for example by knowing a sufficiently nice generating set, presentation or more generally a partial resolution.

\paragraph{This article}
The aim of this article is to give an overview of progress that has been made since Church--Farb--Putman's article \cite{CFP:stabilityconjectureunstable} appeared.
On the one hand, this text collects vanishing and non-vanishing results for the high-dimensional cohomology of various duality groups.
On the other hand, it describes partial resolutions of the dualising modules that were used to obtain some of these results.

\cref{sec_sln}, \cref{sec_mcg} and \cref{sec_aut_F_n} contain the original conjectures of Church--Farb--Putman and explain their current status.
\cref{sec_sl_ring}, \cref{sec_sp} and \cref{sec_chevalley} give generalisations to different types of Chevalley groups that have recently emerged.
\cref{sec_finite_index} summarises results for further groups such as congruence subgroups and concludes with an overview.

\paragraph{Related work}
This article is intended as a comparably short overview that gives pointers to further literature but does not contain many definitions or proofs. 
A more accessible account of some of the results presented here is given in lecture notes \cite{Patzt} by Patzt--Wilson. 
These include background material, exercises and open problems.

\section{$\SL{n}{\mbZ}$}
\label{sec_sln}
For the integral special linear group, \cite[Conjecture 2]{CFP:stabilityconjectureunstable} says the following:

\begin{conjecture}[Church--Farb--Putman]
  \label{conj_cfp}
  $	H^{{n \choose 2}-i}(\SL{n}{\mbZ}; \mbQ) = 0 \text{ for } n \geq i+2$.
\end{conjecture}
Here, ${n \choose 2} = \vcd(\SL{n}{\mbZ})$, so \cref{conj_cfp} claims that the rational cohomology of $\SL{n}{\mbZ}$ vanishes in ``low codimension $i$''.
This vanishing conjecture is equivalent to a high-dimensional stability statement \cite[Conjecture 1]{CFP:stabilityconjectureunstable}: Church--Farb--Putman described a map $\St(\SL{n}{\mbZ})\to \St(\SL{n+1}{\mbZ})$ that induces the homomorphism in the top row of the following commutative diagram.
\begin{equation}
\label{eq_stability_formulation}
\begin{split}
\xymatrix@R=0cm{
 H_i(\SL{n}{\mbZ}; \St(\SL{n}{\mbZ})\otimes \mbQ) \ar[r]& H_{i}(\SL{n+1}{\mbZ}; \St(\SL{n+1}{\mbZ})\otimes \mbQ) \\
\cong & \cong \\
 H^{{n \choose 2}-i}(\SL{n}{\mbZ}; \mbQ) \ar[r]& H^{{n+1 \choose 2}-i}(\SL{n+1}{\mbZ}; \mbQ)
}
\end{split}
\end{equation}
The vertical isomorphisms are a case of Borel--Serre Duality as stated in \cref{eq:BE_duality} and define the homomorphism in the bottom row.
Using the definitions of all these maps, one can see that
\cref{conj_cfp} holds if the horizontal ``stabilisation'' maps in \cref{eq_stability_formulation} are isomorphisms for $n\geq i+2$.\footnote{This holds because the square of this map, so the concatenation $H^{{n \choose 2}-i}(\SL{n}{\mbZ}; \mbQ)\to H^{{n+1 \choose 2}-i}(\SL{n+1}{\mbZ}; \mbQ) \to H^{{n+2 \choose 2}-i}(\SL{n+2}{\mbZ}; \mbQ)$, is the zero map.}

\subsection{Results}
\begin{figure}[t]
\includegraphics{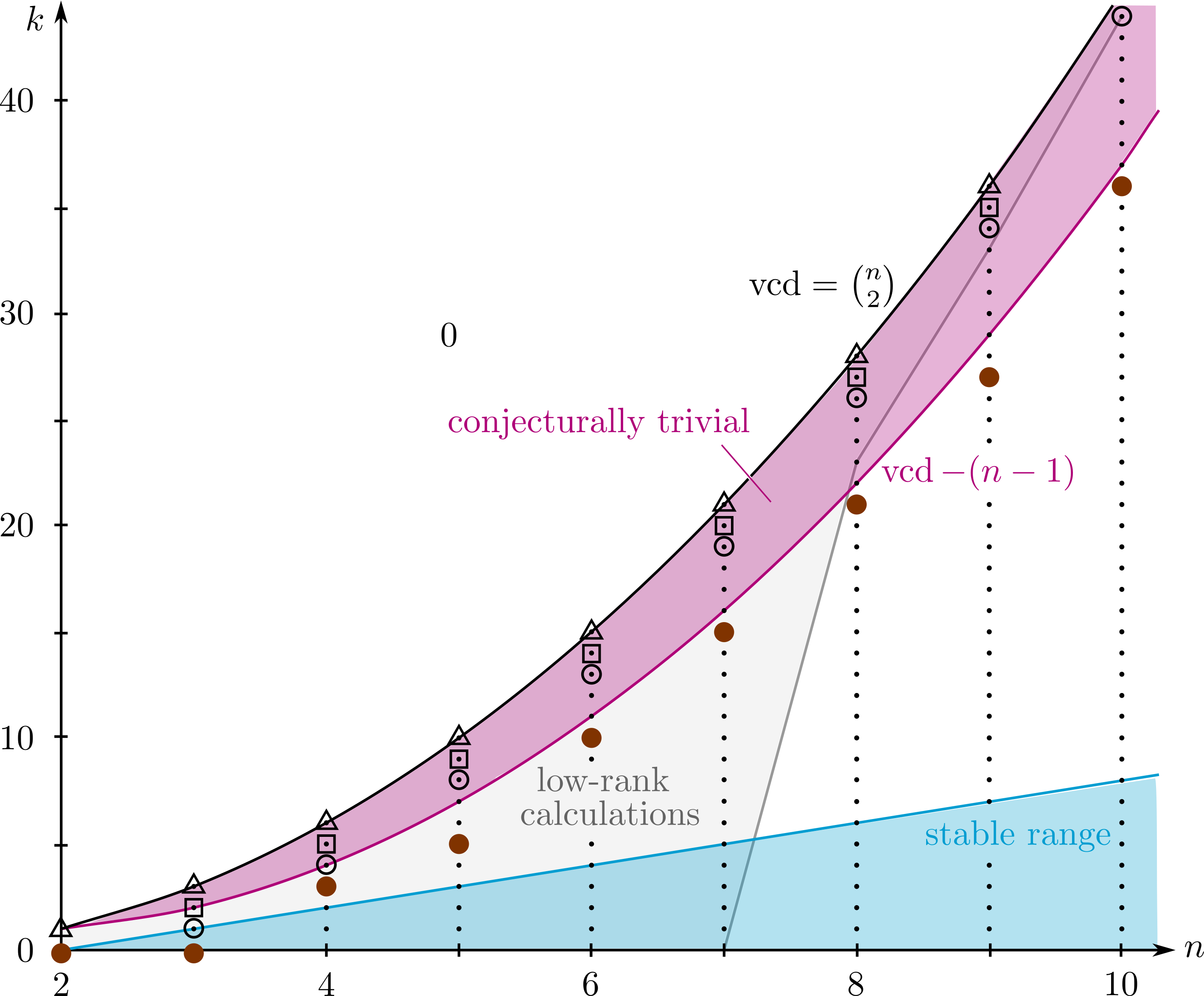}
\caption{$H^k(\SL{n}{\mbZ};\mbQ)$. 
The pink area ${n \choose 2} - (n-2)\leq k \leq {n\choose 2}$ is the range of \cref{conj_cfp}.
Triangle, square and circle mark existing vanishing results (\cref{eq_lee_szcarba}, \cref{eq_church_putman} and \cref{eq_bpmsw});
brown the known non-trivial classes of highest degree (\cref{eq_brown} for $n\geq 8$).}\label{fig:sln_overview}
\end{figure}

\paragraph{Low rank}
For small $n$, the cohomology groups $H^k(\SL{n}{\mbZ};\mbQ)$ can be explicitly calculated using computers. At the moment, this has been done for $n\leq 7$ and for some small $k$ when $n\leq 10$ \cite{Soule1978,Lee1978,ElbazVincent2013, Sikiric2019}.
These low-rank computations are all consistent with \cref{conj_cfp}.

\paragraph{Low codimension}
That \cref{conj_cfp} holds for codimension $i=0$ was proven by Lee--Szczarba \cite{Lee1976} who more generally showed that if $\ring$ is a Euclidean number ring, then the top-degree cohomology of $\SL{n}{\ring}$ is trivial,
\begin{equation}
\label{eq_lee_szcarba}
	H^{\vcd(\SL{n}{\ring})}(\SL{n}{\ring};\mbQ) = 0 \text{ for } \ring \text{ Euclidean and } n\geq 2.
\end{equation}
Church--Farb--Putman \cite{Church2019} reproved Lee--Szczarba's result\footnote{In fact, they proved a stronger statement, see \cite[Theorem C]{Church2019}.} using a ``geometric'' argument that studies certain simplicial complexes (see \cref{sec_resolutions_Steinberg}).
With similar techniques, Church--Putman \cite{CP:codimensiononecohomology} subsequently showed that \cref{conj_cfp} also holds for codimension $i=1$,
\begin{equation}
\label{eq_church_putman}
	H^{{n \choose 2}-1}(\SL{n}{\mbZ};\mbQ) = 0 \text{ for } n\geq 3.
\end{equation}
The currently strongest known bound is due to Br\"uck--Miller--Patzt--Sroka--Wilson \cite{Brueck2022} who showed the case $i=2$ of \cref{conj_cfp}, with a slightly stronger bound on $n$,
\begin{equation}
\label{eq_bpmsw}
	H^{{n \choose 2}-2}(\SL{n}{\mbZ};\mbQ) = 0 \text{ for } n\geq 3.
\end{equation}

\paragraph{Limitations for the vanishing}
Combining recent results by Ash \cite{Ash2024} and Brown \cite{Brown2023}, one can see that if \cref{conj_cfp} holds, then it is sharp for $n \equiv 0,1,2,4 \mod 6$ and close to being sharp otherwise. Their results show that
\begin{align*}
H^{{n \choose 2}-(n-1)}(\SL{n}{\mbZ};\mbQ) &\neq 0 \text{ for } n \equiv 0,1,2,4 \mod 6, \\
H^{{n \choose 2}-n}(\SL{n}{\mbZ};\mbQ) &\neq 0 \text{ for } n \equiv 3,5 \mod 6.
\end{align*}
Low-rank computations show that $H^{{3 \choose 2}-2}(\SL{3}{\mbZ};\mbQ)$ and $H^{{5 \choose 2}-4}(\SL{5}{\mbZ};\mbQ)$ are trivial, so the classes from \cref{eq_brown} are in the highest possible degree for $n\leq 7$ \cite[Table 5]{Sikiric2019}, see also \cref{fig:sln_overview}.
However, Brown--Chan--Galatius--Payne \cite[Corollary 1.9]{Brown2024a} proved that $\dim_{\mbQ}H^{n^2 - (n-1)-k}(\SL{n}{\mbZ};\mbQ)$ grows at least exponentially  with $n$ for all but finitely many $k\in \mbN$. In particular,
\begin{equation}
\label{eq_brown}
	H^{{n \choose 2}-(n-1)-k}(\SL{n}{\mbZ};\mbQ) \neq 0 \text{ for } k\in \mbN\sm\ls 2,3,4,7,8,11,12,16,20,24 \rs,\, n\gg 0.
\end{equation} 
Setting $k=0$, this shows that the bound in \cref{conj_cfp} is asymptotically sharp.

\subsection{Techniques: Partial resolutions of the Steinberg module}
\label{sec_resolutions_Steinberg}
The vanishing results in \cref{eq_lee_szcarba}, \cref{eq_church_putman} and \cref{eq_bpmsw}
were all obtained by applying Borel--Serre Duality (\cref{eq:BE_duality}) and then showing that $H_i(\SL{n}{\mbZ}; \St(\SL{n}{\mbZ})\otimes \mbQ)$ is trivial for appropriate $i$ and $n$. This was done using specific partial resolutions of the Steinberg module $\St(\SL{n}{\mbZ})$.

\subsubsection{How a resolution helps to show vanishing}
Assume that we have a partial flat resolution of $\St(\SL{n}{\mbZ})$, i.e.~an exact sequence
\begin{equation*}
C_k \to C_{k-1}\to \cdots \to C_1 \to C_0 \to \St(\SL{n}{\mbZ}) \to 0,
\end{equation*}
where each $C_i$ is a flat $\SL{n}{\mbZ}$-module. 
Then it follows from standard facts in group homology \cite[Lemma 3.1]{CP:codimensiononecohomology} that $H_*(\SL{n}{\mbZ}; \St(\SL{n}{\mbZ})\otimes \mbQ)$ is given by the homology of the chain complex
\begin{equation*}
C_k \otimes_{\SL{n}{\mbZ}} \mbQ \to C_{k-1}\otimes_{\SL{n}{\mbZ}} \mbQ \to \cdots \to C_1\otimes_{\SL{n}{\mbZ}} \mbQ \to C_0\otimes_{\SL{n}{\mbZ}} \mbQ\to 0.
\end{equation*}
In particular, if for some $i$ and $n$, the coinvariants $C_i \otimes_{\SL{n}{\mbZ}} \mbQ$ vanish, then one gets
\begin{equation*}
	0 = H_i(\SL{n}{\mbZ}; \St(\SL{n}{\mbZ})\otimes \mbQ) \cong H^{{n \choose 2}-i}(\SL{n}{\mbZ};\mbQ).
\end{equation*}
This yields a sufficient criterion for showing \cref{conj_cfp}, which was used to prove it for codimensions $i=0,1,2$.

\subsubsection{What the resolutions look like}
\label{sec_resolutions_form}
In all the partial resolutions considered here, $C_i$ is a free abelian group that is finitely generated as an $\SL{n}{\mbZ}$-module. It has an explicit description of the finitely many orbits of generators and each such generator has finite stabiliser (this is related to flatness \cite[Lemma 3.2]{CP:codimensiononecohomology}).
A partial resolution of length 0 is just a surjective equivariant map $C_0\twoheadrightarrow \St(\SL{n}{\mbZ})$. This is essentially the same as a generating set of $\St(\SL{n}{\mbZ})$ that is invariant under the action of $\SL{n}{\mbZ}$ and has only finitely many orbits of generators. Similarly, a partial resolution of length 1 gives a presentation of $\St(\SL{n}{\mbZ})$. We will describe the partial resolutions that were used for proving \cref{conj_cfp} in codimensions $i=0,1$ in this form.

\paragraph{Buildings and apartments}
For later use, we first consider the more general situation where $\ring$ is a Dedekind domain with fraction field $\field = \on{frac}(R)$ and $\chevalley$ a Chevalley--Demazure group scheme. We then specialise to $\chevalley = \on{SL}_n$, $\ring = \mbZ$, $\field = \mbQ$.
The Steinberg module $\St(\chevalley(\ring))$ is the top-degree homology of the spherical building $\building(\chevalley(\ring))$ associated to $\chevalley(\field)$.
This building is a highly symmetric simplicial complex equipped with an action of $\chevalley(\field)$. It can be described as a union of certain subcomplexes, its apartments. Each apartment $\Sigma$ is isomorphic to the Coxeter complex of the (finite) Weyl group of $\chevalley$ and hence its geometric realisation is a sphere of dimension $\rk(\chevalley)-1$. Its fundamental class gives an element $[\Sigma] \in \St(\chevalley(\ring)) = \tildee{H}_{\rk(\chevalley)-1}(\building(\chevalley(\ring))$.
By the Solomon--Tits Theorem (\cite{Sol:Steinbergcharacterfinite}, \cite[Theorem 4.127]{AB:Buildings}) , these classes generate $\St(\chevalley(\ring))$. 

\begin{figure}
\centering
\includegraphics{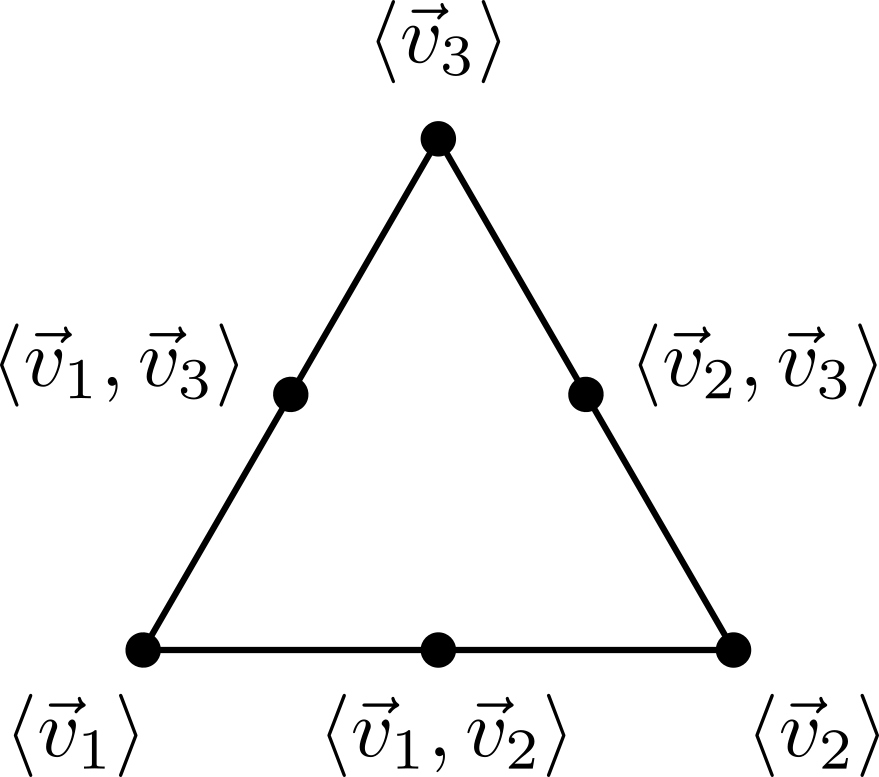}
\caption{An apartment in $\building(\SL{3}{\ring})$, corresponding to a basis $\vec v_1, \vec v_2, \vec v_3$ of $\field^{3}$.}
\label{fig_apartment}
\end{figure}
\paragraph{$\boldsymbol{\on{SL}_n}$}
The building $\building(\SL{n}{\ring})$ is the order complex of the poset of non-zero proper subspaces of $\field^n$. An apartment is determined by a tuple\footnote{The equivalence class of such a tuple is sometimes called a (higher dimensional) ``modular symbol'' \cite{AR:modularsymbolcontinued}.}
$[v_1, \dots, v_n]$, where each $v_i = \langle \vec v_i \rangle$ is a line such that $\vec v_1, \dots, \vec v_n$ is a basis of $\field^{n}$.\footnote{The apartment of  $[v_1, \dots, v_n]$ is the full subcomplex on subspaces spanned by proper subsets of $\ls\vec v_1, \dots, \vec v_n\rs$, see \cref{fig_apartment}. This is isomorphic to the barycentric subdivision of the boundary of an $(n-1)$-simplex, the Coxeter complex of type $\tA_{n-1}$.}
If $\ring$ is Euclidean, not all apartments are needed to generate $\St(\SL{n}{\ring})$: Ash--Rudolph \cite{AR:modularsymbolcontinued} showed that already the classes of \emph{integral apartments}, those $[v_1, \dots, v_n]$ where
$\vec v_1, \dots, \vec v_n$ can be chosen as a basis of $\ring^{n}$, form a generating set. 
Taking $C_0$ to be the free abelian group on all (oriented) integral apartment classes gives a length-0 partial resolution $C_0\to \St(\SL{n}{\ring})$.
Because $\SL{n}{\ring}$ acts transitively on the set of integral apartments, $C_0$ is cyclic as an $\SL{n}{\ring}$-module. This description allows one to show that the coinvariants $C_0\otimes \mbQ$ are trivial, which implies vanishing of the top-degree cohomology of $\SL{n}{\ring}$ in \cref{eq_lee_szcarba}.

\begin{figure}
\centering
\includegraphics{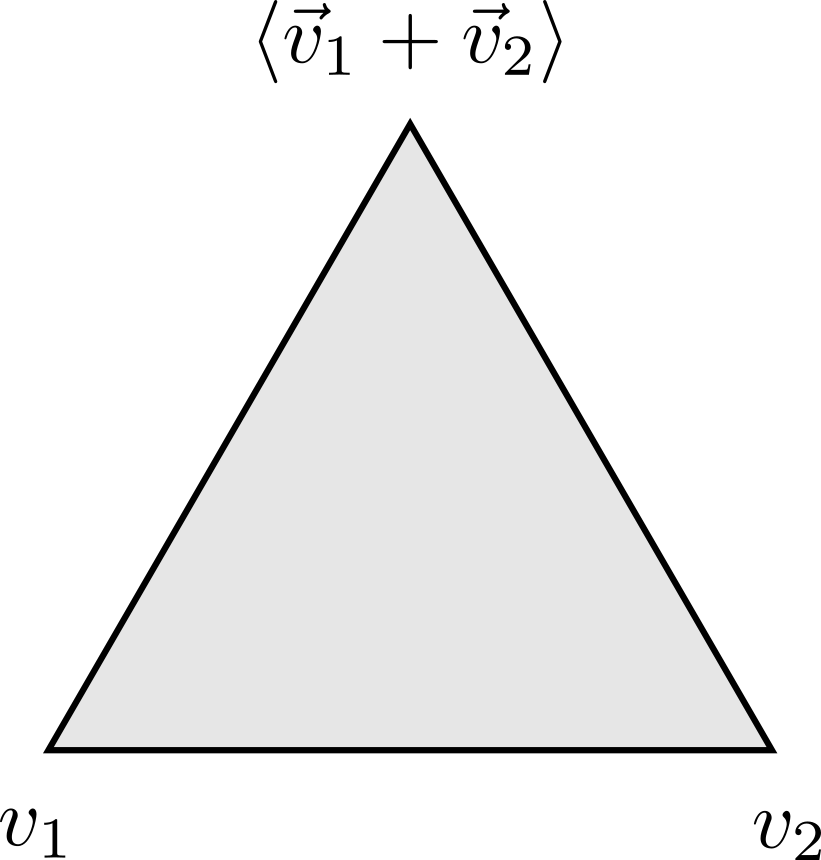}
\caption{A 2-simplex in $\on{BA}_2$ whose boundary leads to Relation 2.~in $\St(\SL{2}{\mbZ})$. Its 1-faces correspond to apartment classes $[v_1, v_2]$, $[v_1, \langle \vec v_1+\vec v_{2} \rangle]$, $[\langle \vec v_1+\vec v_{2} \rangle, v_2]$ and are contained in $\on{B}_2$.}
\label{fig_relation_SL}
\end{figure}

Church--Farb--Putman \cite{Church2019} reproved this integral generation result in a geometric way, that is via showing that a certain simplicial complex $\on{B}_n$ is $(n-2)$-connected. The $(n-1)$-simplices of $\on{B}_n$ are in 1-to-1 correspondence with the integral apartments of $\building(\SL{n}{\ring})$ and $C_0$ is given by the module of simplicial $(n-1)$-chains.
Using similar techniques, Church--Putman \cite[Theorem B]{CP:codimensiononecohomology} also reproved a result by Bykovski\u{\i} \cite{Bykovskiui2003} that gives all relations between integral apartment classes. This yields a presentation, i.e.~length-1 partial resolution $C_1 \to C_0\to \St(\SL{n}{\mbZ})$.
It describes $\St(\SL{n}{\mbZ})$ as the free abelian group generated by symbols $[v_1, \ldots, v_n]$ modulo the following relations:
\begin{equation}
\label{eq_presentation_sl}
\begin{split}
1.\,  [v_1, \ldots, v_n] &= 	(-1)^{\on{len}(\pi)}\cdot 	[v_{\pi(1)},\dots, v_{\pi(n)}] \,\forall \pi \in \Sym_n;\\
2.\, [v_1, \ldots, v_n] &= [v_1, \langle\vec v_1 +\vec 	v_{2} \rangle, v_3, \dots, v_n]+[\langle \vec v_1+\vec v_{2} \rangle, 	v_{2}, v_3, \dots, v_n].
\end{split}
\end{equation}
Both relations in \eqref{eq_presentation_sl} show up as the boundaries of simplices in a simplicial complex $\on{BA}_n\supseteq \on{B}_n$, see \cref{fig_relation_SL}. Showing that $\on{BA}_n$ is $(n-1)$-connected is the main technical difficulty in Church--Putman's proof for the above presentation. 

To prove the codimension $i=2$ case of \cref{conj_cfp} (see \cref{eq_bpmsw}), Br\"uck--Miller--Patzt--Sroka--Wilson  extended this presentation to a length-$2$ partial resolution of $\St(\SL{n}{\mbZ})$ \cite[Section 1.2]{Brueck2022}. 
For this, they showed that a certain simplicial complex $\on{BAA}_n\supseteq \on{BA}_n$ is $n$-connected. To handle the combinatorial complexity of  $\on{BAA}_n$, their argument involves computer calculations.

\section{$\MCG(\Sigma_g)$}
\label{sec_mcg}
The rational cohomology of the surface mapping class group $\on{MCG}(\Sigma_g)$ agrees with that of the moduli space $\mcM_g$ of genus-$g$ Riemann surfaces.
Church--Farb--Putman's analogue of \cref{conj_cfp} in this setup \cite[Conjecture 9]{CFP:stabilityconjectureunstable} asked whether
\begin{equation}
\label{eq_conjecture_MCG}
	H^{4g-5-i}(\on{MCG}(\Sigma_g);\mbQ) \cong H^{4g-5-i}(\mcM_g;\mbQ) = 0 \text{ for } g\gg i ?
\end{equation}
This conjecture, which turned out to be false in general, would have also followed from a conjecture by Kontsevich, see \cite[Remark 7.5]{Morita2015}.

Harer \cite{Har:virtualcohomologicaldimension} and Ivanov \cite{Iva:ComplexescurvesTeichmullerb} showed that the group $\on{MCG}(\Sigma_g)$ is a virtual duality group of dimension $\vcd(\on{MCG}(\Sigma_g))= 4g-5$. More precisely, there are isomorphisms 
\begin{equation}
\label{eq:MCG_duality}
H^{4g-5-i}(\on{MCG}(\Sigma_g);\mbQ) \cong H_i(\on{MCG}(\Sigma_g); D_g \otimes \mbQ),
\end{equation}
where the dualising module $D_g$ is given by the only non-trivial homology group of the curve complex.
Similarly to the case of $\SL{n}{\mbZ}$ above, this gives access to the high-dimensional cohomology of $\on{MCG}(\Sigma_g)$ if one can find appropriate partial resolutions of the dualising module $D_g$.
Broaddus \cite{Bro:Homologycurvecomplex} found a resolution of $D_g$ and gave an explicit generating set for this module. Church--Farb--Putman \cite{CFP:rationalcohomologymapping} used this generating set to show that \cref{eq_conjecture_MCG} is true in codimension $i=0$, which was also shown independently by Morita--Sakasai--Suzuki \cite{Morita2013},
\begin{equation}
\label{eq_vanishing_MCG}
	H^{4g-5}(\on{MCG}(\Sigma_g);\mbQ) \cong H^{4g-5}(\mcM_g;\mbQ) = 0 \text{ for } g \geq 2.
\end{equation}
However, Chan--Galatius--Payne \cite{Chan2021} and later Payne--Willwacher \cite{Payne2021, Payne2024} found many non-trivial classes in the high-dimensional cohomology of $\mcM_g$, showing that
\begin{equation*}
	H^{4g-5-i}(\on{MCG}(\Sigma_g);\mbQ) \neq 0 \text{ for } i\in \ls 1,\ldots,54 \rs \sm\ls 2,5,8,21,52 \rs ,\, g \gg 0.
\end{equation*}
In particular, \cref{eq_conjecture_MCG} is false for all these codimensions $i$.

Since then, related techniques also yielded high-dimensional cohomology classes for  mapping class groups of marked surfaces \cite{Chan2019} and of handlebodies \cite{Hainaut2023,Borinsky2023a, Petersen2024}.

\section{Automorphisms of free groups}
\label{sec_aut_F_n}
Among the three families of groups considered in \cite{CFP:stabilityconjectureunstable}, the situation is least clear for the group $\Aut(F_n)$ of automorphisms of the free group and the related group $\Out(F_n)$ of outer automorphisms.
These groups have virtual cohomological dimensions $\vcd(\Aut(F_n)) = 2n-2$ and $\vcd(\Out(F_n)) = 2n-3$.
Morita \cite[page 390]{Morita1999} conjectured a top-degree vanishing phenomenon for $\Out(F_n)$,
\begin{equation}
\label{eq_conj_Morita}
	H^{2n-3}(\Out(F_n);\mbQ) = 0 \text{ for } n\geq 2 ?
\end{equation}
Church--Farb--Putman did not conjecture vanishing but rather that the high-dimensional cohomology would only depend on the parity of $n$ \cite[Conjecture 12]{CFP:stabilityconjectureunstable},
\begin{equation}
\label{eq_conj_CFP_Aut}
	H^{2n-2-i}(\Aut(F_n);\mbQ), \, H^{2n-3-i}(\Out(F_n);\mbQ) \text{ determined by } n+2\mbZ \text{ for } n\gg i ?
\end{equation}

After these conjectures were made, Bartholdi \cite{Bar:rationalhomologyouter} computed the rational cohomology of $\Out(F_7)$ and showed that it is non-trivial in the top degree,
\begin{equation*}
	H^{11}(\Out(F_7);\mbQ) \cong \mbQ.
\end{equation*}
This disproved Morita's conjecture in \cref{eq_conj_Morita}. 
While \cref{eq_conj_CFP_Aut} might still be true, Bartholdi's result makes this seem much less likely: It implies that a high-dimensional stabilisation for $\Out(F_n)$ as claimed by \cref{eq_conj_CFP_Aut} can only occur for $n\geq 7$. In contrast to that, the analogous stabilisation for $\SL{n}{\mbZ}$ and $\MCG(\Sigma_g)$ already holds for $n\geq 2$ (\cref{eq_lee_szcarba}) and $g\geq 2$ (\cref{eq_vanishing_MCG}).

Bestvina--Feighn \cite{BF:topologyinfinity} showed that $\Aut(F_n)$ and $\Out(F_n)$ are virtual duality groups.
This indicates that it might be possible to study their high-dimensional cohomology with similar techniques as for $\SL{n}{\mbZ}$ and $\MCG(\Sigma_g)$ described above.
The problem with this is that until today, no sufficiently explicit descriptions of the dualising modules are known. 
Hatcher--Vogtmann \cite{HV:complexfreefactors} suggested that the top-degree homology of the free factor complex could be the dualising module of $\Aut(F_n)$, but Himes--Miller--Nariman--Putman \cite{Himes2022} showed that
this is not the case, at least for $n=5$. 
Further candidates for geometric descriptions of the dualising modules are studied in \cite{BG:Homotopytypecomplex} and \cite{BSV:bordificationouterspace}.

\section{$\SL{n}{\ring}$}
\label{sec_sl_ring}
In the past years, vanishing patterns as in \cref{conj_cfp} have been studied for several further families of groups ``similar to $\SL{n}{\mbZ}$''. 
One such family is $\SL{n}{\ring}$, where $\ring$ is a number ring different from $\mbZ$.

\subsection{Vanishing}
Lee--Szczarba's vanishing result for the top-degree cohomology of $\SL{n}{\ring}$ is already formulated in such a setting and applies to all Euclidean number rings $\ring$ (see \cref{eq_lee_szcarba}).

Kupers--Miller--Patzt--Wilson \cite{Kupers2022} proved that this vanishing extends to codimension $i=1$ for two specific such Euclidean number rings, namely the Gaussian integers $\numberring_{-1}$ and the Eisenstein integers $\numberring_{-3}$:
\begin{equation}
\label{eq:Gaussian_Eistenstein}
	H^{(n^2-n)-1}(\SL{n}{\ring};\mbQ) =0 \text{ for } \ring \in \ls \numberring_{-1}, \numberring_{-3} \rs \text{ and } n\geq 3,
\end{equation}
where $n^2-n = \vcd(\SL{n}{\ring})$. To show this, they proved that for these rings, the Steinberg module $\St(\SL{n}{\ring})$ has a ``generalised Bykovski\u{\i} presentation'' that looks like the presentation of $\St(\SL{n}{\mbZ})$ in \cref{eq_presentation_sl}. 
They also showed that this presentation does not hold if $\ring$ is a quadratic number ring that is Euclidean but not additively generated by units \cite[Theorem C]{Kupers2022}.

\subsection{Non-vanishing}
If one drops the assumption that $\ring$ be Euclidean, there is no vanishing in general, not even in top degree.
Recall that the class group $\class(\ring)$ is the abelian group given by isomorphism classes of rank-1 projective $\ring$-modules with the tensor product. A ring $\ring$ is a principal ideal domain (PID) if its class group is trivial, $|\class(\ring)| = 1$.

\paragraph*{$\boldsymbol{\on{SL}_n}$ over non-PIDs}
Church--Farb--Putman \cite{Church2019} showed that $\SL{n}{\ring}$ has non-trivial  top-degree cohomology if $\ring$ is a number ring that is not a PID. More precisely, they proved that 
\begin{equation}
\label{eq:non_vanishing_SLn}
	\dim_\mbQ H^{\vcd(\SL{n}{\ring})}(\SL{n}{\ring}; \mbQ) \geq (|\class(\ring)|-1)^{(n-1)} \text{ for } n \geq 2.
\end{equation}
Their proof also shows that if $\ring$ is a Dedekind domain with $2\leq|\class(\ring)|<\infty$, then $\St(\SL{n}{\ring})$ is not generated by integral apartment classes for $n\geq 2$ \cite[Theorem B]{Church2019}.

\begin{figure}[t]
\includegraphics[width=.6\textwidth]{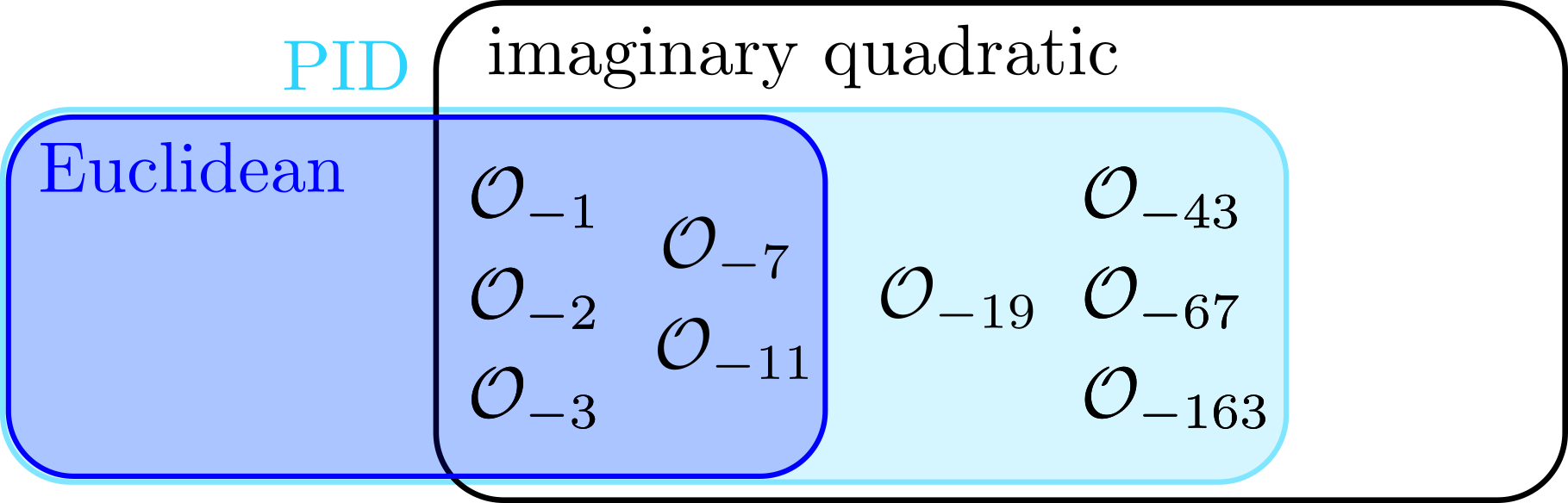}
\caption{Number rings under the assumption of GRH.}\label{fig_number_rings}
\end{figure}

\paragraph*{$\boldsymbol{\on{SL}_n}$ over non-Euclidean PIDs} Every Euclidean ring is a PID. Hence, combining the vanishing result of Lee--Szczarba (\cref{eq_lee_szcarba}) with \cref{eq:non_vanishing_SLn}, the only number rings for which it is unknown whether $\SL{n}{\ring}$ has vanishing top-degree cohomology are non-Euclidean PIDs. Assuming the Generalised Riemann Hypothesis (GRH), there are only four such number rings \cite{Weinberger1973}, see \cref{fig_number_rings}. 
These are the rings of integers $\numberring_{d}$ in $\mbQ(\sqrt{d})$ for  $d\in\ls -19, -43, -67 ,-163 \rs$. Miller--Patzt--Wilson--Yasaki \cite{MPWY:Nonintegralitysome} showed that for the last three, the top-degree cohomology of $\on{SL}_{2n}$ does not vanish,
\begin{equation}
\label{eq_Miller_Yasaki}
	H^{\vcd(\SL{n}{\numberring_{d}})}(\SL{n}{\numberring_{d}}; \mbQ) \neq 0 \text{ for } n \text{ even and }d \in \ls -43, -67, -163 \rs.
\end{equation}
They also refined Church--Farb--Putman's non-integrality result from the previous paragraph in the case of number rings: They showed that if $\ring$ is a number ring and GRH holds, then $\St(\SL{n}{\ring})$ is generated by integral apartments if and only if $\ring$ is Euclidean \cite[Corollary 1.2]{MPWY:Nonintegralitysome}.

\section{$\Sp{2n}{\ring}$}
\label{sec_sp}
Another way to find groups similar to $\SL{n}{\mbZ}$ is to replace $\on{SL}_n$ by other group schemes. Here, in particular the family of symplectic groups $\Sp{2n}{\ring}$ has been studied.

\subsection{Vanishing}
It follows from work of Gunnells that if $\ring$ is a Euclidean number ring, then the cohomology of $\Sp{2n}{\ring}$ vanishes in its virtual cohomological dimension,
\begin{equation}
\label{eq_brueck_sroka}
	H^{\vcd(\Sp{2n}{\ring})}(\Sp{2n}{\ring}; \mbQ) = 0 \text{ for } \ring \text{ Euclidean and } n\geq 1.
\end{equation}
This analogue of Lee--Szczarba's $\on{SL}_n$-result (\cref{eq_lee_szcarba}) first appeared in work of Br\"uck--Santos Rego--Sroka \cite{Brueck2022c}, but was supposedly known to experts before. 
For the case $\ring = \mbZ$, Br\"uck--Patzt--Sroka \cite{Brueck2023, Sroka2021} gave an alternative proof that is independent of Gunnells' work (see \cref{sec_res_symplectic}) and rather studies certain simplicial complexes. They later \cite{Brueck2023a} used these techniques 
to show that cohomology vanishing also occurs in codimension $1$, 
\begin{equation}
\label{eq_brueck_patzt_sroka}
H^{n^2-1}(\Sp{2n}{\mbZ}; \mbQ) = 0 \text{ for } n\geq 2.
\end{equation}

\subsection{Partial resolutions of $\St(\Sp{2n}{\ring})$}
\label{sec_res_symplectic}
Similarly to the setting of $\on{SL}_n$ (see \cref{sec_resolutions_Steinberg}), the vanishing results in \cref{eq_brueck_sroka} and \cref{eq_brueck_patzt_sroka} were obtained using Borel--Serre Duality, which here says that
\begin{equation*}
	 H^{\vcd(\Sp{2n}{\ring})-i}(\Sp{2n}{\ring};\mbQ)  \cong H_i(\Sp{2n}{\ring}; \St(\Sp{2n}{\ring})\otimes \mbQ),
\end{equation*}
together with explicit partial resolutions of the Steinberg module $\St(\Sp{2n}{\ring})$. 

To describe these partial resolutions, let $\ring$ again be a Dedekind domain and $\field$ its fraction field. 
Let $\omega$ be a symplectic form on $\field^{2n}$, i.e.~a non-degenerate alternating bilinear form $\omega: \field^{2n}\times \field^{2n} \to \field$. 
The symplectic group $\Sp{2n}{\ring}$ is given by all elements of $\on{GL}_{2n}(\ring)\subseteq \on{GL}_{2n}(\field)$ that preserve $\omega$.
A subspace $V\subseteq \field^{2n}$ is called isotropic if $\omega|_{V\times V}$ is trivial. A basis $\vec v_1, \vec v_{\bar{1}}, \dots, \vec v_n,\vec v_{\bar{n}}$ of $\field^{2n}$ is called symplectic if
\begin{equation*}
	\omega(\vec v_i,\vec v_j) = \omega(\vec v_{\bar{i}},\vec v_{\bar{j}}) = 0\quad\text{and}\quad \omega(\vec v_i, \vec v_{\bar{j}}) = \delta_{ij},
\end{equation*} 
where $\delta_{ij}$ is the Kronecker delta.

The Steinberg module $\St(\Sp{2n}{\ring})$ is the $(n-1)$-st homology group of the building $\building(\Sp{2n}{\ring})$.
This building is the poset of non-zero isotropic subspaces of $\field^{2n}$.
An apartment in $\building(\Sp{2n}{\ring})$ is determined by a tuple $[v_1, \dots, v_{\bar{n}}]$, where each $v_i = \langle \vec v_i \rangle$ is a line such that 
$\vec v_1, \dots, \vec v_{\bar{n}}$ is a symplectic basis of $\field^{2n}$.\footnote{The apartment of  $[v_1, \dots, v_{\bar{n}}]$ is the full subcomplex on all isotropic subspaces spanned by subsets of $\vec v_1, \dots, \vec v_{\bar{n}}$. This is isomorphic to the barycentric subdivision of the boundary of an $n$-dimensional cross polytope, the Coxeter complex of type $\tB_{n}$.}
It is \emph{integral} if $\vec v_1, \dots, \vec v_{\bar{n}}$ is a basis of $\ring^{2n}$. Gunnells \cite{Gun:Symplecticmodularsymbols} showed that $\St(\Sp{2n}{\ring})$ is generated by integral apartments for $\ring$ Euclidean. Br\"uck--Sroka \cite{Brueck2023} gave a new proof for $\ring = \mbZ$. With Patzt \cite{Brueck2023a}, they extended the techniques to obtain a presentation of $\St(\Sp{2n}{\mbZ})$. 
They showed \cite[Theorem B]{Brueck2023a} that it is the free abelian group generated 
by symbols $[v_1, \dots, v_{\bar{n}}]$ 
modulo the following relations:
\begin{equation}
\label{eq_presentation_sp}
\begin{split}
\hspace*{-0.05cm}1 .\, [v_1, \dots, v_{\bar{n}}] = &\,
	(-1)^{\on{len}(\pi)}\cdot 
	[v_{\pi(1)},\dots, v_{\pi(\bar{n})}] \,\forall \pi \in \CoxeterGroupTypeB_n;\\
\hspace*{-0.05cm}2.\, [v_1, \dots, v_{\bar{n}}] = &\, [v_1, \langle\vec v_1 +\vec 
	v_{\bar{1}} \rangle, v_2, \dots, v_{\bar{n}}]+[\langle \vec v_1+\vec v_{\bar{1}} \rangle, 
	v_{\bar{1}},v_2, \dots,  v_{\bar{n}}]; \\
\hspace*{-0.05cm}3.\, [v_1, \dots, v_{\bar{n}}] = &\, [v_1,\langle \vec v_{\bar{1}} - 
	\vec v_{\bar{2}} \rangle, \langle \vec v_1 + \vec v_2 \rangle, v_{\bar{2}}, v_3, \dots, 
	v_{\bar{n}}] \\
	 &\, + [\langle \vec v_{\bar{1}}-\vec 
	v_{\bar{2}} \rangle , v_2, \langle\vec v_1 + \vec v_2\rangle, v_{\bar{1}}, v_3, \dots, 
	 v_{\bar{n}}].
\end{split}
\end{equation}

\begin{figure}
\centering
\includegraphics{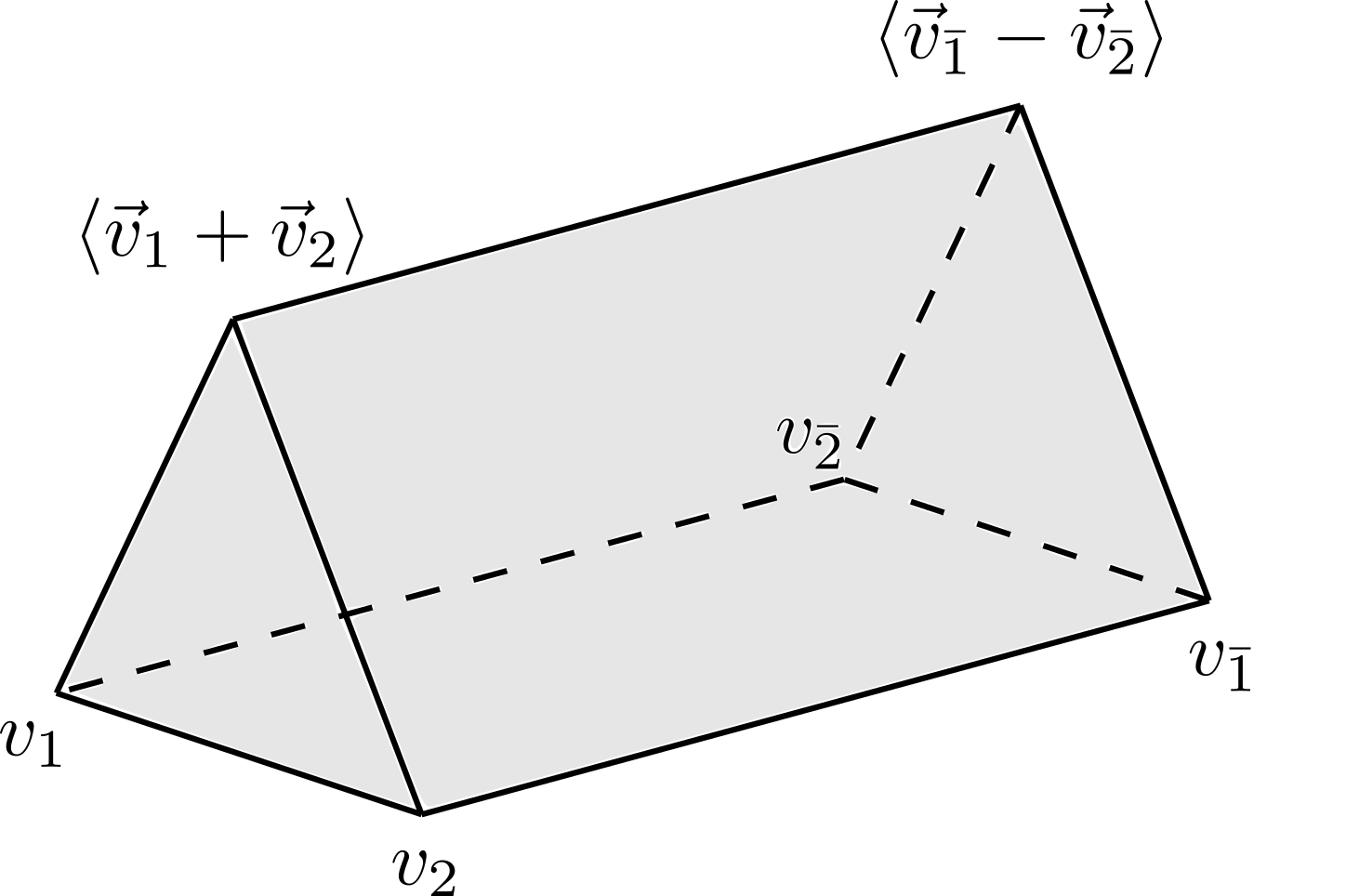}
\caption{A subcomplex of $\on{IAA}_2$ whose boundary leads to Relation 3.~in $\St(\Sp{4}{\mbZ})$.}
\label{fig_relation_Sp}
\end{figure}

\noindent
Here, $\CoxeterGroupTypeB_n$ is the signed permutation group, the Coxeter group of type $\tB_n = \tC_n$.
Similar to the case of $\SL{n}{\mbZ}$, these three relations occur as boundaries in an $n$-connected complex $\on{IAA}_n$, see \cref{fig_relation_Sp}.

\subsection{Non-vanishing}
Different sources of non-trivial cohomology classes provide limitations for a high-dimensional vanishing of $H^*(\Sp{2n}{\ring};\mbQ)$.

\paragraph{$\boldsymbol{\on{Sp}_{2n}(\mbZ)}$ and the moduli space $\mcA_n$}
The rational cohomology of $\Sp{2n}{\mbZ}$ equals that of $\mcA_n$, the moduli space of principally polarised abelian varieties of dimension $2n$. A result by van der Geer \cite{Geer1999} implies that
\begin{equation}
\label{eq_brandt_etal}
H^{n^2 - n}(\Sp{2n}{\mbZ};\mbQ) \cong H^{n^2-n}(\mcA_n;\mbQ) \neq 0 \text{ for }n\geq 1.
\end{equation}
In fact, the dimension $\dim_{\mbQ}H^{n^2 - n - k}(\Sp{2n}{\mbZ};\mbQ)$ even grows at least exponentially with $n$ for $k=0$ and all but finitely many $k\geq 1$ \cite[Corollary 1.6]{Brown2024a}.

\paragraph{$\on{Sp}_{2n}$ over non-PIDs}
Br\"uck--Himes \cite{Brueck2023b} gave an analogue of Church--Farb--Putman's non-vanishing result in \cref{eq:non_vanishing_SLn}. 
They proved that for number rings $\ring$, the top-degree cohomology of $\Sp{2n}{\ring}$ is non-trivial if $\ring$ is not a PID,
\begin{equation}
\label{eq_non_trivial_cohomology_sp}
	\dim_\mbQ H^{\vcd(\Sp{2n}{\ring})}(\Sp{2n}{\ring}; \mbQ) \geq (|\class(\ring)|-1)^n \text{ for } n \geq 1.
\end{equation}
They also showed that if $\ring$ is a Dedekind domain with $2\leq |\class(\ring)|<\infty$, then $\St(\Sp{2n}{\ring})$ is not generated by integral apartment classes for $n\geq 1$ \cite[Theorem 1.2]{Brueck2023b}.

\section{Chevalley groups}
\label{sec_chevalley}
For the top degree, i.e.~the cohomology in the virtual cohomological dimension, there is a vanishing result that covers a much bigger class of groups than just $\on{SL_n}$ and $\on{Sp}_{2n}$.

Let $\ring$ be a number ring and $\chevalley$ a Chevalley--Demazure group scheme. 
As mentioned in \cref{sec_resolutions_form}, the Steinberg module $\St(\chevalley(\ring))$ is generated by the fundamental classes of apartments in the building $\building(\chevalley(\ring))$. Call such a class \emph{integral} if it is a $\chevalley(\ring)$-translate of the class of the standard apartment\footnote{This designated apartment in $\building(\chevalley(\ring))$ is uniquely determined by the choice of a BN-pair in $\chevalley(\field)$, see \cite[Chapter 6]{AB:Buildings}.}. This notion generalises integrality of apartment classes in the special linear and symplectic settings described above.

T\'oth \cite{Toth2005} showed that if $\ring$ is Euclidean, then for almost all Chevalley types, the Steinberg module $\St(\chevalley(\ring))$ is generated by integral apartment classes.
Br\"uck--Santos Rego--Sroka \cite{Brueck2022c} deduced from this that the top-degree cohomology of $\chevalley(\ring)$ is trivial,
\begin{equation}
\label{eq_cohomology_vanishing_specific_types}
	H^{\vcd(\chevalley(\ring))}(\chevalley(\ring); \mbQ ) = 0 \text{ for } R \text{ Euclidean, }\chevalley \text{ type } \tA_n,\, \tB_n,\, \tC_n,\, \tD_n,\, \tE_6,\,\tE_7.
\end{equation}

They then stated the following conjecture as \cite[Question 1.2]{Brueck2022c}.

\begin{conjecture}
\label{conj_brueck}
Let $\ring$ be a Euclidean number ring and $\chevalley$ a Chevalley-Demazure group scheme. Then
\begin{equation*}
	H^{\vcd(\chevalley(\ring))-i}(\chevalley(\ring);\mbQ) = 0 \text{ for } i< \rk(\chevalley).
\end{equation*}
\end{conjecture}

For $\chevalley = \on{SL}_n$ and $\ring = \mbZ$, this specialises to \cref{conj_cfp} by Church--Farb--Putman.

\paragraph*{Known cases}
For $i=0$, \cref{conj_brueck} is known to hold for all irreducible $\chevalley$ except those of type $\tF$, $\tG$ or $\tE_8$ (\cref{eq_cohomology_vanishing_specific_types}).
For $i=1$, it it is presently known to be true for four families: for $\SL{n}{\mbZ}$ (\cref{eq_church_putman}); for $\SL{n}{\numberring_{-1}}$ and $\SL{n}{\numberring_{-3}}$ (\cref{eq:Gaussian_Eistenstein}); and for $\Sp{2n}{\mbZ}$ (\cref{eq_brueck_patzt_sroka}).
For $i=2$, the only known instance is $\SL{n}{\mbZ}$ (\cref{eq_bpmsw}).
\cref{conj_brueck} is consistent with low-rank cohomology computations for $\SL{n}{\mbZ}$, $n\leq 10$ \cite{Soule1978,Lee1978,ElbazVincent2013, Sikiric2019}; for $\Sp{2n}{\mbZ}$, $n\leq 4$ \cite{Igusa1962, Hain2002, Hulek2012}; and for $\SL{n}{\ring}$, $n\leq 3$, where $\ring$ is a Euclidean ring of integers in any imaginary quadratic field \cite{Schwermer1983, DutourSikiric2016}.

\paragraph*{Limitations}
\cref{eq:non_vanishing_SLn} and \cref{eq_non_trivial_cohomology_sp} show that the vanishing claimed by \cref{conj_brueck} can only occur for rings $\ring$ that are PIDs. 
As mentioned above, $\ring$ being Euclidean only excludes four further rings under the assumption of the Generalised Riemann Hypothesis. This stronger condition seems natural given the non-vanishing results of Miller--Patzt--Wilson--Yasaki (\cref{eq_Miller_Yasaki}).

Known cohomology classes show that one also cannot hope for vanishing results in higher codimensions:  The highest degree classes known for $\SL{n}{\mbZ}$ (\cref{eq_brown}) and for $\Sp{2n}{\mbZ}$ 
(\cref{eq_brandt_etal}) are exactly one below the conjectured vanishing range, in codimensions $(n-1) = \rk(\on{SL}_n)$ and $n = \rk(\on{Sp}_{2n})$, respectively.

\section{Further results and overview}

\paragraph{Congruence subgroups}
\label{sec_finite_index}
The high-dimensional vanishing results presented above depend on the existence of torsion elements in the groups under questions. These are used to show vanishing of the coinvariants of partial resolutions as described in \cref{sec_resolutions_Steinberg} (see e.g.~\cite[Proof of Theorem 3.4]{Brueck2022c} or \cite[Section 3.2, p.~1012]{CP:codimensiononecohomology}).
In particular, vanishing is not preserved under passing to finite index subgroups. 
For example, the top-degree cohomology of congruence subgroups in $\SL{n}{\mbZ}$ is non-trivial and its rank grows quickly with $n$ \cite{Paraschivescu1997,MPP:topdimensionalcohomology,Schwermer2021, Scalamandre2023}.
The situation is similar for congruence subgroups in $\MCG(\Sigma_g)$ \cite{Fullarton2020, BBP:highdimensionalcohomology}.
Church--Farb--Putman conjectured that there might be some kind of high-dimensional stability for congruence subgroups in $\SL{n}{\mbZ}$  nonetheless.
In \cite{Miller2020}, Miller--Nagpal--Patzt stated a precise version of this conjecture \cite[Conjecture 6.1]{Miller2020} using representation stability and proved it in some instances.

\paragraph{$\on{GL}_n$ and configuration spaces}
For vanishing and non-vanishing results in the top-degree cohomology of $\GL{n}{\ring}$, see \cite{Putman2022a}.
High-dimensional stability patterns in the context of configuration spaces are studied in \cite{Maguire2016,Knudsen2023,Yameen2023}.

\paragraph{Overview}
\label{sec_overview}
\begin{table}
\renewcommand{\arraystretch}{1.3}
\centering
\begin{tabular}{c|c|c|c}
Group $\group$ & $i=0$& $i=1$ & $i=2$\\
\hline
$\SL{n}{\mbZ}$ & Yes \cite{Lee1976} & Yes \cite{CP:codimensiononecohomology} & {Yes \cite{Brueck2022}}\\
\hline
$\on{MCG}(\Sigma_g)$  & \makecell{{Yes} {\cite{CFP:rationalcohomologymapping}} } & \makecell{{No} \cite{Chan2021}} & \\
\hline
$\on{Out}(F_n)$ & \makecell{{Probably no  \cite{Bar:rationalhomologyouter}}\\{(low rank calculations)}}  & &\\
\hline
\hline
\makecell{$\SL{n}{\ring}$ \\$\ring$ Eucl.} &   {Yes} \cite{Lee1976} & \makecell{Yes for\\ $ \ring \in \ls \numberring_{-1}, \numberring_{-3} \rs$ \cite{Kupers2022}}& \\ 
\hline
\makecell{$\SL{n}{\ring}$ \\$\ring$ not Eucl.} &   \makecell{No if $R\neq \numberring_{-19}$ \cite{Church2019, MPWY:Nonintegralitysome} \\ (assuming GRH)}&&   \\
\hline
$\Sp{2n}{\mbZ}$ &  \makecell{{Yes} {\cite{Brueck2023}}}&   \makecell{{Yes} \cite{Brueck2023a}}& \\
\hline
\makecell{$\Sp{2n}{\ring}$ \\$\ring$ not PID } & \makecell{{ No  \cite{Brueck2023b}}} &&  \\
\hline
\makecell{$\chevalley_n(\ring)$ \\ $\ring$ Eucl.} &  \makecell{Yes for type $\tA_n$, $\tB_n$, \\ {$\tC_n$, $\tD_n$, $\tE_6$, $\tE_7$} \cite{Brueck2022c}}&&\\
\hline
\hline
$\on{MCG}(V_g)$ &  & No \cite{Hainaut2023,Borinsky2023a} &\\
\hline
\makecell{congruence \\ subgroups } & \makecell{No \\ \cite{Paraschivescu1997,MPP:topdimensionalcohomology,Schwermer2021, Scalamandre2023,Fullarton2020, BBP:highdimensionalcohomology}} &&
\end{tabular}
\caption{A table showing whether stably, $H^{\vcd(\group)-i}(\group;\mbQ) = 0$ for different groups $\group$. $\ring$ denotes a number ring, $\Sigma_g$ and $V_g$ a genus-$g$ surface and handlebody, $F_n$ the free group of rank $n$.}
\label{tab_overview_CFP_conjecture}
\end{table}

Of the original three conjectures in \cite{CFP:stabilityconjectureunstable}, only the one for $\SL{n}{\mbZ}$ still has a good chance of being true.
However, results for further arithmetic groups indicate that there might be a vanishing phenomenon similar to that of \cref{conj_cfp} that occurs in a larger class of groups.
The precise form of this still has to be determined, one candidate is given by \cref{conj_brueck}.
The high-dimensional (non-)vanishing results mentioned in this article are summarised in \cref{tab_overview_CFP_conjecture}.


\begin{ack}
I would like to thank Jeremy Miller, Jennifer C.~H.~Wilson and an anonymous referee for their comments, which helped to improve the content and presentation of this article.
\end{ack}

\begin{funding}
This work was partially supported by the Deutsche Forschungsgemeinschaft (DFG, German Research Foundation) -- Project-ID 427320536 -- SFB 1442, as well as by Germany’s Excellence Strategy EXC 2044 -- 390685587, Mathematics Münster: Dynamics–Geometry–Structure.
\end{funding}

\bibliographystyle{ems}
\bibliography{mybibliography}









\end{document}